# WEAK CONVERGENCE OF MEASURE-VALUED PROCESSES AND *R*-POINT FUNCTIONS

By Mark Holmes[1] and Edwin Perkins[2]

*Technische Universiteit Eindhoven and University of British Columbia*

We prove a sufficient set of conditions for a sequence of finite measures on the space of cadlag measure-valued paths to converge to the canonical measure of super-Brownian motion in the sense of convergence of *finite-dimensional distributions*. The conditions are convergence of the Fourier transform of the *r*-point functions and perhaps convergence of the "survival probabilities." These conditions have recently been shown to hold for a variety of statistical mechanical models, including critical oriented percolation, the critical contact process and lattice trees at criticality, all above their respective critical dimensions.

**1. Motivation.** In the last few years, a number of rescaled models from interacting particle systems and statistical physics have been shown to converge to the canonical measure of super-Brownian motion. The models include critical oriented percolation above four dimensions [6], critical contact processes above four dimensions [5] and critical lattice trees above eight dimensions [7], all for sufficiently spread-out kernels. In each of these cases, what is actually proved is convergence of the Fourier transforms of the moment measures (or *r*-point functions). Our modest objective here is to translate this result into the more conventional probabilistic language of weak convergence of stochastic processes. To those well versed in weak convergence arguments, we fear this may be one of the proverbial much-needed gaps in the literature, but to others who have complained to us, it is an irritant that should be spelled out once and for all.

Received January 2006.
[1]Supported in part by a UGF from the University of British Columbia and by the Netherlands Organization for Scientific Research.
[2]Supported in part by an NSERC research grant.
*AMS 2000 subject classifications.* Primary 60G57, 60K35; secondary 60F05.
*Key words and phrases.* *r*-point functions, measure-valued processes, super-Brownian motion, canonical measure, critical oriented percolation.

The limiting measure is a sigma-finite measure (not a probability) on the space of continuous measure-valued paths, which presents some additional minor worries. The full convergence on path space remains open in all of the above settings due to the absence of any tightness result on path space. Even the natural statement of convergence of finite-dimensional distributions requires convergence of the survival probabilities (see Proposition 2.4 below), a result which was only recently discovered for critical oriented percolation [3, 4] and which is currently being pursued in the other contexts mentioned above. So, in the end, we thought that someone should advertise this state of affairs and we have acquiesced in the writing of this note. If you are reading this in a journal, at least one editor and/or referee has agreed with us.

**2. Introduction.** Consider a discrete-time, critical nearest-neighbor branching random walk on $\mathbb{Z}^d$, starting with a single particle at the origin. That is, at time $n \in \mathbb{Z}_+$, each individual gives birth to a random number of offspring, each of which immediately takes a step to a randomly chosen nearest neighbor of its parent. Assume that each parent dies immediately after giving birth, that the offspring distribution has mean one and finite variance $\gamma > 0$, and that each of the offspring laws and random walk steps are independently chosen.

Extend the branching random walk to all times $t \geq 0$ by making it a right-continuous step function. Let $\mathcal{M}_t = \{x_t^{(\alpha)} : \alpha \in I_t\}$ denote the set of locations of particles in $\mathbb{Z}^d$ which are alive at time $t$. We have suppressed the details of the labeling system (see, e.g., Section II.3 in [8]), but as multiple occupancies are allowed, some labeling scheme is needed here.

In order to describe the scaling limit, we represent the model as a *cadlag* measure-valued process by setting

$$X_t^n = \frac{C_1}{n} \sum_{\alpha \in I_{nt}} \delta_{x_{nt}^{(\alpha)}/(C_2\sqrt{n})},$$

where $C_1 = \gamma^{-1/2}$ and $C_2 = d^{-1/2}$. If $E$ and $M$ are separable metric spaces, then $M_F(E)$ denotes the space of finite Borel measures on $E$ with the topology of weak convergence and $D(M)$ denotes the space of *cadlag* $M$-valued paths with the Skorokhod topology. With probability 1, $X_t^n$ is a finite measure for all $n \in \mathbb{Z}_+$ and $t \geq 0$, so that $\{X_t^n\}_{t \geq 0} \in D \equiv D(M_F(\mathbb{R}^d))$.

The *extinction time* $S : D \to [0, \infty]$ is defined by

$$S(X) \equiv \inf\{s > 0 : X_s = 0_M\},$$

where $0_M$ is the zero measure on $\mathbb{R}^d$ and $\inf \varnothing = \infty$. Next, we define a sequence of measures $\mu_n \in M_F(D)$ by

$$(1) \qquad \mu_n(\bullet) \equiv C_3 n \mathbb{P}(\{X_t^n\}_{t \geq 0} \in \bullet),$$

where $C_3 = 1$ for this branching random walk model.

Let $M_\sigma(D)$ denote the $\sigma$-finite measures on $D$ which assign finite mass to $\{S > \varepsilon\}$ for all $\varepsilon > 0$, with the topology of weak convergence defined as follows.

DEFINITION 2.1 (*Weak convergence*). Let $\{\nu_n : n \in \mathbb{N} \cup \{\infty\}\} \subset M_\sigma(D)$. We write $\nu_n \stackrel{w}{\Longrightarrow} \nu_\infty$ if for every $\varepsilon > 0$,

$$\nu_n^\varepsilon(\bullet) \equiv \nu_n(\bullet, S > \varepsilon) \stackrel{w}{\Longrightarrow} \nu_\infty(\bullet, S > \varepsilon) \equiv \nu_\infty^\varepsilon(\bullet) \qquad \text{as } n \to \infty,$$

where the convergence is in $M_F(D)$.

It is a standard result in the superprocess literature (see, e.g., [8], Theorem II.7.3) that there exists $\mathbb{N}_0 \in M_\sigma(D)$, supported by the continuous paths in $D$ which remain at $0_M$ after time $S$, and called the *canonical measure of super-Brownian motion* (CSBM), such that $\mu_n \stackrel{w}{\Longrightarrow} \mathbb{N}_0$. In [8], one is working with branching Brownian motion instead of branching random walk but it is trivial to modify the arguments. We have chosen our constants $C_i$ so that the branching and diffusion parameters of our limiting super-Brownian motion are both equal to one. Much is known about $\mathbb{N}_0$, for example, as in Theorem II.7.2(iii) of [8], we have for every $b > 0$ that

$$(2) \qquad \mathbb{N}_0(X_b(1) \in A \setminus \{0\}) = \left(\frac{2}{b}\right)^2 \int_A e^{-(2/b)x}\, dx.$$

Let $l \geq 1$ and $\vec{t} = \{t_1, \ldots, t_l\} \in [0, \infty)^l$. We use $\pi_{\vec{t}} : D \to M_F(\mathbb{R}^d)^l$ to denote the projection map satisfying $\pi_{\vec{t}}(X) = (X_{t_1}, \ldots, X_{t_l})$. The *finite-dimensional distributions* of $\nu \in M_\sigma(D)$ are the measures $\nu^\varepsilon \pi_{\vec{t}}^{-1} \in M_F(M_F(\mathbb{R}^d)^l)$ given by

$$\nu^\varepsilon \pi_{\vec{t}}^{-1}(H) \equiv \nu^\varepsilon(\{X : \pi_{\vec{t}}(X) \in H\}), \qquad H \in \mathcal{B}(M_F(\mathbb{R}^d)^l).$$

DEFINITION 2.2 (*Convergence of f.d.d.*). Let $\{\nu_n : n \in \mathbb{N} \cup \{\infty\}\} \subset M_\sigma(D)$. We write $\nu_n \stackrel{f.d.d.}{\Longrightarrow} \nu_\infty$ if for every $\varepsilon > 0$, $m \in \mathbb{N}$ and $\vec{t} \in [0, \infty)^m$,

$$\nu_n^\varepsilon \pi_{\vec{t}}^{-1} \stackrel{w}{\Longrightarrow} \nu_\infty^\varepsilon \pi_{\vec{t}}^{-1} \qquad \text{as } n \to \infty,$$

where the convergence is in $M_F(M_F(\mathbb{R}^d)^m)$.

If $\nu_\infty$ is supported by continuous paths in $D$, it is easy to see that weak convergence to $\nu_\infty$ (Definition 2.1) implies convergence of the finite-dimensional distributions to $\nu_\infty$ (Definition 2.2). An additional tightness condition is needed for the converse.

We now work in a more abstract setting than the branching random walk described above, in which $\{\mu_n\}$ is any sequence of finite measures on $D$. For $k \in \mathbb{R}^d$, let $\phi_k(x) = e^{ik \cdot x}$ and write $E_{\mu_n}[Y]$ for $\int Y \, d\mu_n$ and $X_t(\phi)$ for $\int \phi(x) X_t(dx)$, respectively. Consider the following convergence condition on the moment measures of $\mu_n$:

$$E_{\mu_n}\left[\prod_{i=1}^{r-1} X_{t_i}(\phi_{k_i})\right] \to E_{\mathbb{N}_0}\left[\prod_{i=1}^{r-1} X_{t_i}(\phi_{k_i})\right]$$

(3)
$$\text{for } r \geq 2, \vec{t} \in [0, \infty)^{r-1}, \vec{k} \in \mathbb{R}^{d(r-1)}.$$

An explicit formula for the right-hand side of (3) can be found in Section 1.2.3 of [6].

Of course, (3) does hold for the $\mu_n$ defined in (1) for branching random walk, but our interest in this condition arises from a number of models in which $\mathcal{M}_t$ is the (finite) set of occupied sites in $\mathbb{Z}^d$ at time $t$. Examples include the critical contact process, critical oriented percolation or critically weighted lattice trees, all with the natural definitions of "occupied site." For $r \geq 2$ and $\vec{t} \in [0, \infty)^{r-1}$, the $r$-point functions for this model are $B_{\vec{t}}(\vec{x}) = \mathbb{P}(x_i \in \mathcal{M}_{t_i}, i = 1, \ldots, r-1)$, while the $\hat{r}$-*point functions* are the Fourier transforms of these quantities,

$$\widehat{B}_{\vec{t}}(\vec{k}) = \sum_{\vec{x}} e^{i\vec{k} \cdot \vec{x}} B_{\vec{t}}(\vec{x}), \qquad \vec{k} \in \mathbb{R}^{d(r-1)},$$

which are defined whenever $B_{\vec{t}}(\vec{x})$ is summable in $\vec{x}$. Define $X_t^n \in M_F(\mathbb{R}^d)$ by

$$X_t^n \equiv \frac{C_1}{n} \sum_{x: C_2\sqrt{n}x \in \mathcal{M}_{nt}} \delta_x$$

and assume that $\mu_n$ given by (1) defines a finite measure on $D$. An easy calculation then shows that

$$\frac{C_1^{r-1} C_3}{n^{r-2}} \widehat{B}_{n\vec{t}}\left(\frac{\vec{k}}{C_2\sqrt{n}}\right) = E_{\mu_n}\left[\prod_{i=1}^{r-1} X_{t_i}(\phi_{k_i})\right]$$

whenever $B_{n\vec{t}}(\vec{x})$ is summable. Therefore, the asymptotic formulae for the $\hat{r}$-point functions for sufficiently spread-out critical rescaled oriented percolation ($d > 4$), critical rescaled lattice trees ($d > 8$), and critical rescaled contact processes ($d > 4$) derived in [6, 7] and work in progress in [5], respectively, immediately imply (3) in each of these cases. Moreover, in each of these models, it is known that $\mu_n$ is a finite measure supported by $D$, as is required above.

In what follows, we use $\mathcal{D}_F$ to denote the set of discontinuities of a function $F$. A function $Q: M_F(\mathbb{R}^d)^m \to \mathbb{R}$ is called a *multinomial* if $Q(\vec{X})$ is a real multinomial in $\{X_1(1), \ldots, X_m(1)\}$. A function $F: M_F(\mathbb{R}^d)^m \to \mathbb{C}$ is said to be *bounded by a multinomial* ($|F| \leq Q$) if there exists a multinomial $Q$ such that $|F(\vec{X})| \leq Q(\vec{X})$ for every $\vec{X} \in M_F(\mathbb{R}^d)^m$. The main results of this paper are the following two propositions. By the above, the first result is applicable in each of the three settings [5, 6, 7].

PROPOSITION 2.3. *Let $\{\mu_n\}_{n \geq 1}$ be a sequence of finite measures on $D(M_F(\mathbb{R}^d))$ such that* (3) *holds. Then for every $s > 0$, $\lambda > 0$, $m \geq 1$, $\vec{t} \in [0, \infty)^m$ and every Borel measurable $F: M_F(\mathbb{R}^d)^m \to \mathbb{C}$ bounded by a multinomial and such that $\mathbb{N}_0 \pi_{\vec{t}}^{-1}(\mathcal{D}_F) = 0$, we have*

$$E_{\mu_n}[X_s(1)F(\vec{X}_{\vec{t}})] \to E_{\mathbb{N}_0}[X_s(1)F(\vec{X}_{\vec{t}})] \tag{4}$$

*and*

$$E_{\mu_n}[F(\vec{X}_{\vec{t}})I_{\{X_s(1) > \lambda\}}] \to E_{\mathbb{N}_0}[F(\vec{X}_{\vec{t}})I_{\{X_s(1) > \lambda\}}]. \tag{5}$$

For critical oriented percolation above the critical spatial dimension of four (and for sufficiently spread-out kernels), [3, 4] show that

$$\mu_n(S > \varepsilon) \to \mathbb{N}_0(S > \varepsilon) \quad \text{for every } \varepsilon > 0. \tag{6}$$

The corresponding results for critical lattice trees and critical contact processes are conjectured to be true above the critical dimension; the latter is currently work in progress (see [5] for the contact process). The next result allows us to strengthen the conclusion of Proposition 2.3 under (6).

PROPOSITION 2.4. *Let $\{\mu_n\}_{n \geq 0}$ be a sequence of finite measures on $D(M_F(\mathbb{R}^d))$ such that* (3) *and* (6) *hold. Then $\mu_n \overset{f.d.d.}{\Longrightarrow} \mathbb{N}_0$.*

In particular, the results of [3, 4, 6], together with Proposition 2.4, imply that above the critical dimension and at the critical occupation probability, the scaling limit (in the sense of finite-dimensional distributions) of sufficiently spread-out oriented percolation is CSBM. Tightness, and hence a full statement of weak convergence, remains an open problem.

The additional condition (6) is necessary (consider the test function 1) because $\mu_n$ and $\mathbb{N}_0$ are unnormalized. In [1], a conditional limit theorem for rescaled lattice trees above eight dimensions is proved in which the limit distribution (ISE) is $\mathbb{N}_0(\int_0^\infty X_s \, ds \in \cdot \mid \int_0^\infty X_s(1) \, ds = 1)$. The conditioning means that all of the involved measures are probabilities and so (6) is not needed.

The following assumption will be in force throughout the rest of the paper.

ASSUMPTION 2.5. $\mathcal{F}$ denotes a class of $\mathbb{C}$-valued bounded continuous functions that is closed under conjugation, is convergence determining for $M_F(\mathbb{R}^d)$ and contains the constant function 1.

We show in Section 4 that both propositions are consequences of standard exponential moment bounds for $\mathbb{N}_0$ and the following theorem. By convention, an empty product is defined to be 1.

THEOREM 2.6. *Let $\mu_n, \mu \in M_F(D(M_F(\mathbb{R}^d)))$. Suppose that for every $l \in \mathbb{Z}_+$ and $\vec{t} \in [0, \infty)^l$, we have:*

1. *there exists a $\delta = \delta(\vec{t}) > 0$ such that for all $\theta_i < \delta$, $E_\mu[e^{\sum_{i=1}^l \theta_i X_{t_i}(1)}] < \infty$;*

2. *for every $\phi_i \in \mathcal{F}$,*

$$E_{\mu_n}\left[\prod_{i=1}^l X_{t_i}(\phi_i)\right] \to E_\mu\left[\prod_{i=1}^l X_{t_i}(\phi_i)\right] < \infty. \tag{7}$$

*Then for every $l \in \mathbb{N}$ and every $\vec{t} \in [0, \infty)^l$, $\mu_n \pi_{\vec{t}}^{-1} \xrightarrow{w} \mu \pi_{\vec{t}}^{-1}$.*

Note that some of the $t_i$'s may be the same in (7).

The remainder of this paper is organized as follows. In Section 3, we prove Theorem 2.6. In Section 4, we prove Propositions 2.3 and 2.4.

**3. Proof of Theorem 2.6.** In this section, we prove Theorem 2.6 as a consequence of Lemmas 3.2–3.7. Lemma 3.2 is standard and states that if a sequence of finite measures is tight, then every subsequence has a further subsequence that converges. Lemma 3.3 establishes tightness of the $\{\mu_n \pi_{\vec{t}}^{-1} : n \in \mathbb{N}\}$ for each fixed $\vec{t}$. Thus, every subsequence of the $\mu_n \pi_{\vec{t}}^{-1}$ has a further subsequence that converges. Lemma 3.4 states that any limit point of the $\{\mu_n \pi_{\vec{t}}^{-1} : n \in \mathbb{N}\}$ must have the same moments (7) as $\mu \pi_{\vec{t}}^{-1}$. Lemma 3.5 extends equality of the moments on the right-hand side of (7) for two measures $\mu, \mu'$ to all $\phi_i \geq 0$ bounded and continuous. Lemmas 3.6 and 3.7 together imply that under condition 1 of Theorem 2.6, equality of the moments in Lemma 3.5 implies equality of the underlying finite measures on $M_F(\mathbb{R}^d)^m$. Taken together, they show that since all subsequential limit points have the same moments (7), the limit points all coincide and thus the whole sequence converges to that limit point. Thus, Theorem 2.6 follows immediately from the lemmas proved in this section.

Recall the notion of *tightness* for finite measures.

DEFINITION 3.1. A set of finite Borel measures $F \subset M_F(E)$ on a metric space $E$ is *tight* if $\sup_{\mu \in F} \mu(E) < \infty$ and for every $\eta > 0$ there exists a compact $K \subset E$ such that $\sup_{\mu \in F} \mu(K^c) < \eta$.

LEMMA 3.2. *If $F \subset M_F(E)$ is tight, then every sequence in $F$ has a subsequence which converges in $M_F(E)$ (weak convergence).*

LEMMA 3.3. *Let $\mu_n, \mu \in M_F(D)$. Suppose that $E_{\mu_n}[1] \to E_\mu[1] < \infty$ and that for every $t \in [0, \infty)$ and every $\phi \in \mathcal{F}$,*

(8) $$E_{\mu_n}[X_t(\phi)] \to E_\mu[X_t(\phi)] < \infty.$$

*Then for each $m \in \mathbb{Z}_+$ and every $\vec{t} \in [0, \infty)^m$, the set of measures $\{\mu_n \pi_{\vec{t}}^{-1} : n \in \mathbb{N}\}$ is tight on $M_F(\mathbb{R}^d)^m$.*

PROOF. That $\sup_n \mu_n \pi_{\vec{t}}^{-1}(M_F(\mathbb{R}^d)^m) < \infty$ for every $m, \vec{t}$ is trivial (as is the $m = 0$ case) since $E_{\mu_n}[1] \to E_\mu[1] < \infty$. It remains to prove the existence of the appropriate compact set for $m \geq 1$.

For $m = 1$, let $\varepsilon > 0$ and $t \geq 0$. Define the mean measures $\nu_n, \nu \in M_F(\mathbb{R}^d)$ by $\nu_n(A) = E_{\mu_n}[X_t(A)]$ and $\nu(A) = E_\mu[X_t(A)]$. Then (8) implies that $\nu_n \to \nu$ in $M_F(\mathbb{R}^d)$ and $\sup_n \nu_n(\mathbb{R}^d) \equiv L < \infty$. Choose $M$ such that $L/M < \varepsilon/2$. Then by Markov's inequality,

(9) $$\sup_n \mu_n(X_t(\mathbb{R}^d) > M) \leq \frac{L}{M} < \frac{\varepsilon}{2}.$$

Fix $\eta > 0$. There exists $K_{-1} \subset \mathbb{R}^d$ compact such that $\nu(K_{-1}^c) < \eta/2$. Furthermore, there exists $K_0 \subset \mathbb{R}^d$ compact such that $\nu(\overline{K_0^c}) \leq \nu(K_{-1}^c)$ [e.g., the set $K_0 = \{x : d(x, K_{-1}) \leq 1\}$]. Since $\nu_n \to \nu$ in $M_F(\mathbb{R}^d)$ and $\overline{K_0^c}$ is closed, we have

$$\limsup_n \nu_n(\overline{K_0^c}) \leq \nu(\overline{K_0^c}) < \frac{\eta}{2}.$$

It follows easily that there exists a compact subset $K$ of $\mathbb{R}^d$ such that $\sup_n \nu_n(K^c) < \eta$. Another application of Markov's inequality implies that

$$\sup_n \mu_n(X_t(K^c) > \eta^{1/4}) \leq \eta^{-1/4} \sup_n E_{\mu_n}[X_t(K^c)] < \eta^{3/4}.$$

Choose $\eta^{1/4} = 2^{-j}$. There then exists $K_j \subset \mathbb{R}^d$ compact such that

(10) $$\sup_n \mu_n\left(X_t(K_j^c) > \frac{1}{2^j}\right) \leq \frac{1}{2^{3j}}.$$

Choose $N$ such that $8^{1-N} < \varepsilon/2$ and let

$$\mathbf{K} \equiv \bigcap_{j \geq N} \left\{Y : Y(K_j^c) \leq \frac{1}{2^j}\right\} \cap \{Y : Y(\mathbb{R}^d) \leq M\}.$$

Now, $\mathbf{K}$ is compact in $M_F(\mathbb{R}^d)$ (see, e.g., the proof of Theorem II.4.1 in [8]) and

$$\mathbf{K}^c = \bigcup_{j \geq N} \left\{Y : Y(K_j^c) > \frac{1}{2^j}\right\} \cup \{Y : Y(\mathbb{R}^d) > M\}.$$

Thus, (10) and (9) imply that

$$\sup_n \mu_n(X_t \in \mathbf{K}^c) \leq \sum_{j=N}^{\infty} \frac{1}{2^{3j}} + \frac{\varepsilon}{2} \leq \frac{1}{8^{N-1}} + \frac{\varepsilon}{2} < \varepsilon, \tag{11}$$

which verifies that the $\mu_n \pi_t^{-1}$, $n \geq 1$ are tight.

For $m > 1$ and $\vec{t} \in [0, \infty)^m$, we have from (11) that for each $i \in \{1, \ldots, m\}$, there exists $\mathbf{K}_i \subset M_F(\mathbb{R}^d)$ compact such that $\sup_n \mu_n \pi_{t_i}^{-1}(\mathbf{K}_i^c) < \varepsilon/m$. Let $\mathbf{K} = \mathbf{K}_1 \times \mathbf{K}_2 \times \cdots \times \mathbf{K}_m$. Then $\mathbf{K} \subset M_F(\mathbb{R}^d)^m$ is compact and

$$\sup_n \mu_n \pi_{\vec{t}}^{-1}(\mathbf{K}^c) \leq \sup_n \sum_{i=1}^{m} \mu_n \pi_{t_i}^{-1}(\mathbf{K}_i^c) < \varepsilon,$$

which gives the result. □

LEMMA 3.4. *Suppose that $\mu_n, \mu \in M_F(D)$ satisfy the second hypothesis of Theorem 2.6. Fix $l \geq 0$ and $\vec{t} \in [0, \infty)^l$. If, for a given subsequence $\mu_{n_k}$, we have $\mu_{n_k} \pi_{\vec{t}}^{-1} \xrightarrow{w} \nu$ in $M_F(M_F(\mathbb{R}^d)^l)$, then for each $\vec{m} \in \mathbb{Z}_+^l$ and $\phi_{ij} \in \mathcal{F}$,*

$$E_\nu \left[ \prod_{i=1}^{l} \prod_{j=1}^{m_i} Y_i(\phi_{ij}) \right] = E_\mu \left[ \prod_{i=1}^{l} \prod_{j=1}^{m_i} X_{t_i}(\phi_{ij}) \right]. \tag{12}$$

PROOF. Assume first that $\nu(1) > 0$. Since $\mu_{n_k} \pi_{\vec{t}}^{-1}(1) \to \nu(1)$, by normalization, we may assume that $\mu_{n_k} \pi_{\vec{t}}^{-1}$ are probabilities on $M_F(\mathbb{R}^d)^l$. Let $\vec{m}$ and $\phi_{ij}$ be as in the lemma and set $W = \prod_{i=1}^{l} \prod_{j=1}^{m_i} X_{t_i}(\phi_{ij})$ and $W_1 = \prod_{i=1}^{l} \prod_{j=1}^{m_i} X_{t_i}(1)$. Condition 2 from Theorem 2.6 implies that

$$\sup_k \mu_{n_k}(|W|^2) \leq C_{\vec{\phi}} \sup_k \mu_{n_k}(W_1^2) < \infty$$

(recall that we can repeat $t_i$'s). The assumed weak convergence implies that $\mu_{n_k} W^{-1} \xrightarrow{w} \nu W^{-1}$ as measures in $M_F(\mathbb{C})$. It follows from a standard result in weak convergence (see, e.g., Proposition 2.3 in the Appendix of [2]) that the left-hand side of (12) is equal to $\lim_{k \to \infty} E_{\mu_{n_k}}[W]$. The same is true of the right-hand side by the second hypothesis in Theorem 2.6, where we use a base vector $\vec{t'}$ with appropriately repeated $t_i$'s.

If $\nu(1) = 0$, then $\mu_{n_k}(D) \to 0$ and so if $W$ is as above, we have

$$\left| \int W \, d\mu_{n_k} \right|^2 \leq \mu_{n_k}(D) \int |W|^2 \, d\mu_{n_k} \to 0,$$

where $L^2$ boundedness of $|W|$ follows as above. Therefore, the right-hand side of (12) is 0 by hypothesis (as above) and thus equals the left-hand side. □

LEMMA 3.5. *Suppose that $l \geq 0$, $\vec{m} \in \mathbb{Z}_+^l$ and $\mu, \mu' \in M_F(M_F(\mathbb{R}^d)^l)$. If*

$$\tag{13} E_\mu\left[\prod_{i=1}^l \prod_{j=1}^{m_i} Y_i(\phi_{ij})\right] = E_{\mu'}\left[\prod_{i=1}^l \prod_{j=1}^{m_i} Y_i(\phi_{ij})\right]$$

*holds (and both quantities are finite) for every $\phi_{ij} \in \mathcal{F}$, then (13) holds for all bounded, continuous $\phi_{ij} \geq 0$.*

PROOF. If $l = 0$ or $\sum m_i = 0$, then the conclusion is trivial, so we may assume that $l > 0$ and $\sum m_i > 0$. Since $1 \in \mathcal{F}$, we have $E_\mu[\prod_{i=1}^l \prod_{j=1}^{m_i} Y_i(1)] < \infty$. Let $\phi_{ij} \in \mathcal{F}$ and $\varphi((x_{ij})) = \prod_{i=1}^l \prod_{j=1}^{m_i} \phi_{ij}(x_{ij})$. Applying Fubini's Theorem to (13), using the fact that the $\phi_{ij} \in \mathcal{F}$ are bounded, we have

$$\tag{14} \begin{aligned}\int \cdots \int \varphi((x_{ij})) E_\mu\left[\prod_{i=1}^l \prod_{j=1}^{m_i} Y_i(dx_{ij})\right] \\ = \int \cdots \int \varphi((x_{ij})) E_{\mu'}\left[\prod_{i=1}^l \prod_{j=1}^{m_i} Y_i(dx_{ij})\right].\end{aligned}$$

We claim that for any $r \geq 1$, the set of functions $\mathcal{F}_r \equiv \{\prod_{i=1}^r \phi_i(x_i) : \phi_i \in \mathcal{F}\}$ is a determining class for $M_F(\mathbb{R}^{dr})$. For real-valued functions, this is Proposition 3.4.6 of [2]. The fact that $\mathcal{F}$ is closed under conjugation easily implies that it is a determining class for complex-valued measures. This allows us to apply the proof in [2] to verify the claim.

Therefore, the products of $\phi_{ij}$ in (14) uniquely determine the measure $\nu$ on $\mathbb{R}^{d\sum m_i}$ defined by $\nu(d\vec{x}) = E_\mu[\prod_{i=1}^l \prod_{j=1}^{m_i} Y_i(dx_{ij})]$. Thus, (14) holds for all $\phi_{ij}$ bounded and continuous. Now, apply Fubini's Theorem again to get (13) for all $\phi_{ij}$ bounded and continuous, as required. □

In the following lemma, $\mathcal{B}_b(\mathbb{R}^d, \mathbb{R}_+)$ denotes the bounded, nonnegative, real-valued functions on $\mathbb{R}^d$, and $\overline{D_0}^{bp}$ denotes the bounded pointwise closure of $D_0 \subset \mathcal{B}_b(\mathbb{R}^d, \mathbb{R}_+)$, that is, the smallest set containing $D_0$ that is closed under bounded pointwise convergence.

LEMMA 3.6. *Suppose that $\mu, \mu' \in M_F(M_F(\mathbb{R}^d)^m)$ and assume that $D_0 \subset \mathcal{B}_b(\mathbb{R}^d, \mathbb{R}_+)$ satisfies $\overline{D_0}^{bp} = \mathcal{B}_b(\mathbb{R}^d, \mathbb{R}_+)$. If for all $h_j \in D_0$,*

$$\tag{15} E_\mu[e^{-\sum_{j=1}^m Y_j(h_j)}] = E_{\mu'}[e^{-\sum_{j=1}^m Y_j(h_j)}],$$

*then $\mu = \mu'$.*

PROOF. By dominated convergence, the identity (15) extends to all bounded, nonnegative, Borel measurable $h_j$. The result follows by using a standard monotone class argument (e.g., see Theorem 4.3 in the Appendix of [2]) on

(16) $$\mathcal{H} \equiv \{\Phi \in \mathcal{B}_b(M_F(\mathbb{R}^d)^m, \mathbb{R}) : E_\mu[\Phi(\vec{Y})] = E_{\mu'}[\Phi(\vec{Y})]\}. \qquad \square$$

LEMMA 3.7. *Let $\mu \in M_F(M_F(\mathbb{R}^d)^m)$. Suppose that there exists a $\delta > 0$ such that for all $\theta_i < \delta$,*

(17) $$E_\mu[e^{\sum_{i=1}^m \theta_i Y_i(1)}] < \infty.$$

*Then for every bounded continuous $0 \leq \psi_i$, the quantity $E_\mu[e^{-\sum_{i=1}^m Y_i(\psi_i)}]$ is uniquely determined by the collection of mixed moments*

$$\left\{ E_\mu\left[\prod_{i=1}^m Y_i(h_i)^{n_i}\right] : 0 \leq h_i \leq 1 \text{ is continuous}, i = 1, \ldots, m \right\}.$$

PROOF. Without loss of generality, we may assume that $m > 0$. Given bounded continuous $\psi_i \geq 0$, define $h_i = \psi_i / \|\psi_i\|_\infty \in [0,1]$ (set $h_i = 0$ if $\psi_i \equiv 0$).

For $\operatorname{Re} z_i < \delta$, $i = 1, \ldots, m$, let

$$f(z_1, \ldots, z_m) = E_\mu[e^{\vec{z} \cdot \vec{Y}(\vec{h})}].$$

Use (17), the Taylor expansion for the exponential function and Fubini's Theorem to see that for $\|\vec{z}\|_\infty < \delta$,

$$f(z_1, \ldots, z_m) = \sum_{l=0}^\infty \frac{1}{l!} E_\mu\left[ \sum_{\vec{n} \in \mathbb{Z}_+^m : \sum n_i = l} \frac{l!}{\prod_{i=1}^m n_i!} \prod_{i=1}^m (z_i Y_i(h_i))^{n_i} \right].$$

Hence, the mixed moments of the form

(18) $$E_\mu\left[\prod_{i=1}^m Y_i(h_i)^{n_i}\right], \qquad n_i \in \mathbb{Z}_+,$$

uniquely determine $f(z)$ for $\|\vec{z}\|_\infty < \delta$. A simple application of dominated convergence and (17) allows us to take the differentiate through the expectation and show that for fixed $z_1, \ldots, z_{j-1}, z_{j+1}, \ldots, z_m$ satisfying $\operatorname{Re} z_i < \delta$ for $i \neq j$, $f(z)$ is analytic in $\operatorname{Re} z_j < \delta$ (and not just $|z_j| < \delta$). Now, use induction on $j \leq m$ to see that moments of the form (18) uniquely determine $f(z_1, \ldots, z_m)$ for $\operatorname{Re} z_1, \ldots, \operatorname{Re} z_{j-1} < \delta$, $|z_j| \vee \cdots \vee |z_m| < \delta$. Here, one uses the aforementioned analyticity in $\operatorname{Re} z_j < \delta$ in the induction step. Apply this result at $z_i = -\|\psi_i\|_\infty$ to complete the proof. $\square$

**4. Applications of Theorem 2.6.** In this section, we prove Propositions 2.3 and 2.4, which relate the asymptotic formulae for the $\widehat{r}$-point functions for various spread-out models above their critical dimensions to the convergence to CSBM. Recall that $\phi_k(x) = e^{ik \cdot x}$. In this section, we fix our convergence-determining class of functions for $M_F(\mathbb{R}^d)$ to be

$$\mathcal{F} = \{\phi_k : k \in \mathbb{R}^d\}, \tag{19}$$

which clearly satisfies Assumption 2.5.

The following lemma will be used to verify the exponential moment hypothesis of Theorem 2.6 for $\mathbb{N}_0$. The branching and diffusion parameters for $\mathbb{N}_0$ are taken to be 1. The lemma is well known, but we include a proof for completeness.

LEMMA 4.1. *For every $b \geq 0$, the following hold:*
1. *for every $\lambda > 0$, $\mathbb{N}_0(X_b(1) = \lambda) = 0$;*
2. *for every $l$ and $\vec{t} \in [0, \infty)^l$, there exists a $\delta = \delta(\vec{t}, b) > 0$ such that for $\theta_i < \delta$,*

$$E_{\mathbb{N}_0}[X_b(1) e^{\sum_{i=1}^{l} \theta_i X_{t_i}(1)}] < \infty; \tag{20}$$

3. *for every $m$ and $\vec{t} \in [0, \infty)^m$ and every $\varepsilon > 0$, there exists a $\delta = \delta(\vec{t}, \varepsilon) > 0$ such that for $\theta_i < \delta$,*

$$E_{\mathbb{N}_0}[e^{\sum_{i=1}^{m} \theta_i X_{t_i}(1)} I_{\{S > \varepsilon\}}] < \infty. \tag{21}$$

PROOF. The first assertion is trivial by (2) and the fact that $\mathbb{N}_0(X_0(1) > 0) = 0$.

As above, we may assume $b > 0$ in part 2. The fact that $X_t = 0_M$ for $t \geq S$ $\mathbb{N}_0$-a.e. implies that for each $\eta > 0$,

$$X_b(1) \leq I_{\{S > b\}} C_\eta e^{\eta X_b(1)}, \qquad \mathbb{N}_0\text{-a.e.}$$

Therefore, part 2 will follow from part 3 with $\varepsilon = b = t_{l+1}$ and $m = l + 1$.

For the last claim of the lemma, we abuse our notation and let $E_{X_0}$ also denote expectation for our standard super-Brownian motion starting at $X_0$. Let $\mathcal{G}_t$ denote the canonical filtration generated by the coordinates $X_s$ of our super-Brownian motion for $s \leq t$. If $H : M_F(\mathbb{R}^d) \to [0, \infty)$ is continuous, then for $t \geq s > 0$,

$$E_{\mathbb{N}_0}[H(X_t) | \mathcal{G}_s] = E_{X_s}[H(X_{t-s})], \qquad \mathbb{N}_0\text{-a.e.} \tag{22}$$

This is easily derived, for example, from the convergence of branching random walk to $\mathbb{N}_0$ mentioned in Section 2, the Markov property for branching random walk and the analogous convergence result for super-Brownian motion (e.g., Theorem II.5.2 of [8]).

We may assume, without loss of generality, that $0 < \varepsilon < t_i < t_{i+1}$ for each $i$. Observe from (22) that

$$E_{\mathbb{N}_0}[e^{\sum_{i=1}^m \theta_i X_{t_i}(1)} I_{\{S>\varepsilon\}}]$$
(23)
$$= E_{\mathbb{N}_0}[E_{X_{t_{m-1}}}[e^{\theta_m X_{t_m - t_{m-1}}(1)}] e^{\sum_{i=1}^{m-1} \theta_i X_{t_i}(1)} I_{\{S>\varepsilon\}}]$$
$$\leq E_{\mathbb{N}_0}[e^{2\theta_m X_{t_{m-1}}(1)} e^{\sum_{i=1}^{m-1} \theta_i X_{t_i}(1)} I_{\{S>\varepsilon\}}],$$

where the inequality holds for $\theta_m$ sufficiently small depending on $t_m - t_{m-1}$, by Lemma III.3.6 of [8]. The last line of (23) has no $t_m$ dependence and, proceeding by induction, it is enough to show that for sufficiently small $\theta > 0$,

(24) $$E_{\mathbb{N}_0}[e^{\theta X_{t_1}(1)} I_{\{S>\varepsilon\}}] < \infty.$$

For $\theta > 0$ small enough [depending on $(\varepsilon, t_1)$], as in (23), the left-hand side is

$$E_{\mathbb{N}_0}[E_{\mathbb{N}_0}[e^{\theta X_{t_1}(1)}|\mathcal{G}_\varepsilon] I_{\{S>\varepsilon\}}] \leq E_{\mathbb{N}_0}[e^{2\theta X_\varepsilon(1)} I_{\{S>\varepsilon\}}]$$
(25)
$$\leq E_{\mathbb{N}_0}[e^{2\theta X_\varepsilon(1)} I_{\{X_\varepsilon(1)>0\}}]$$
$$= \left(\frac{2}{\varepsilon}\right)^2 \int_0^\infty e^{2\theta x} e^{-2x/\varepsilon} \, dx,$$

where the last equality holds by (2). The last line of (25) is finite for sufficiently small $\theta > 0$ (depending on $\varepsilon$) and the result follows. $\square$

PROOF OF PROPOSITION 2.3. Define $\mu_{n,s}, \mathbb{N}_{0,s} \in M_F(D(M_F(\mathbb{R}^d)))$ by

(26)
$$\mu_{n,s}(A) = \int_A X_s(1) \, d\mu_n,$$
$$\mathbb{N}_{0,s}(A) = \int_A X_s(1) \, d\mathbb{N}_0.$$

That these measures are finite follows from the fact that for $s > 0$,

(27) $$\mu_{n,s}(D) = E_{\mu_n}[X_s(1)] \to E_{\mathbb{N}_0}[X_s(1)] < \infty.$$

For all $l \geq 0$ and $\vec{k} \in \mathbb{R}^{dl}$,

$$E_{\mu_{n,s}}\left[\prod_{i=1}^l X_{t_i}(\phi_{k_i})\right] = E_{\mu_n}\left[X_s(1) \prod_{i=1}^l X_{t_i}(\phi_{k_i})\right]$$
(28)
$$\to E_{\mathbb{N}_0}\left[X_s(1) \prod_{i=1}^l X_{t_i}(\phi_{k_i})\right]$$
$$= E_{\mathbb{N}_{0,s}}\left[\prod_{i=1}^l X_{t_i}(\phi_{k_i})\right],$$

where, even in the $l = 0$ case, the presence of the factor $X_s(1)$ ensures that the convergence in (28) follows from (3).

By Lemma 4.1, we have that

(29) $$E_{\mathbb{N}_{0,s}}[e^{\sum_{i=1}^m \theta_i X_{t_i}(1)}] < \infty,$$

for $\theta_i > 0$ sufficiently small depending on $\vec{t}$ and $s$. In view of (27), (28) and (29), we may apply Theorem 2.6 to the measures $\mu_{n,s}, \mathbb{N}_{0,s}$ to obtain

$$\mu_{n,s}\pi_{\vec{t}}^{-1} \stackrel{w}{\Longrightarrow} \mathbb{N}_{0,s}\pi_{\vec{t}}^{-1}.$$

Thus, (4) holds for every bounded continuous $F$. The extension to bounded, Borel-measurable $F$ satisfying $\mathbb{N}_{0,s}\pi_{\vec{t}}^{-1}(\mathcal{D}_F) = 0$ is standard. For $F$ as in the theorem we may assume that $F \geq 0$. The extension to $F$ dominated by a multinomial $Q$ is obtained by an easy uniform integrability argument since $\lim_{n\to\infty} E_{\mu_{n,s}}[Q(\vec{X}_{\vec{t}})] = E_{\mathbb{N}_{0,s}}[Q(\vec{X}_{\vec{t}})]$.

To prove the second claim, we define

$$G_s \equiv \begin{cases} 0, & \text{if } X_s(1) = 0, \\ \dfrac{I_{\{X_s(1) > \lambda\}}}{X_s(1)}, & \text{otherwise.} \end{cases}$$

Then $G_s$ is continuous, except when $X_s(1) = \lambda$, and is bounded above by $\frac{1}{\lambda}$. Thus, Lemma 4.1 and (4) show that for $F$ as in Proposition 2.3,

$$E_{\mu_n}[X_s(1)G_s F(\vec{X}_{\vec{t}})] \to E_{\mathbb{N}_0}[X_s(1)G_s F(\vec{X}_{\vec{t}})],$$

that is,

$$E_{\mu_n}[I_{\{X_s(1) > \lambda\}} F(\vec{X}_{\vec{t}})] \to E_{\mathbb{N}_0}[I_{\{X_s(1) > \lambda\}} F(\vec{X}_{\vec{t}})]. \qquad \square$$

PROOF OF PROPOSITION 2.4. We apply Theorem 2.6 to the finite measures $\mu_n^\varepsilon$ and $\mathbb{N}_0^\varepsilon$ defined by

(30) $$\mu_n^\varepsilon(\bullet) = \mu_n(\bullet, S > \varepsilon), \qquad \mathbb{N}_0^\varepsilon(\bullet) = \mathbb{N}_0(\bullet, S > \varepsilon).$$

Fix $l \in \mathbb{Z}_+$ and $\vec{t} \in [0,\infty)^l$. By Lemma 4.1, for $\delta(\vec{t},\varepsilon) > 0$ sufficiently small and for $\theta_i < \delta$,

$$E_{\mathbb{N}_0^\varepsilon}[e^{\sum_{i=1}^l \theta_i X_{t_i}(\mathbb{R}^d)}] < \infty,$$

so that the first condition of Theorem 2.6 is satisfied. The second condition is trivially true if any $t_i = 0$, so we assume that $t_i > 0$ for each $i$.

Let $\eta > 0$. Fix $l \in \mathbb{Z}_+$, $\vec{k} \in \mathbb{R}^{dl}$ and write $F(\vec{X}_{\vec{t}}(\vec{\phi})) \equiv \prod_{i=1}^l X_{t_i}(\phi_{k_i})$. By hypothesis [repeat $t_i$'s in (3)], we have

$$E_{\mu_n}[F^2(\vec{X}_{\vec{t}}(\vec{1}))] \to E_{\mathbb{N}_0}[F^2(\vec{X}_{\vec{t}}(\vec{1}))] < \infty,$$

so there exists $C_0(\vec{t})$ such that $\sup_n E_{\mu_n}[F^2(\vec{X}_{\vec{t}}(\vec{1}))]^{1/2} \leq C_0$. Choose $\lambda_0 = \lambda_0(\eta, C_0, \varepsilon)$ sufficiently small so that

(31) $$\mathbb{N}_0(X_\varepsilon(1) \in (0, \lambda_0]) < \left(\frac{\eta}{6C_0}\right)^2.$$

By part 2 of Proposition 2.3 with $F \equiv 1$, we have

$$\mu_n(X_\varepsilon(1) > \lambda_0) \to \mathbb{N}_0(X_\varepsilon(1) > \lambda_0).$$

Combining this with (6) gives $\mu_n(X_\varepsilon(1) \in (0, \lambda_0]) \to \mathbb{N}_0(X_\varepsilon(1) \in (0, \lambda_0])$. It follows from (31) that there exists $n_0$ such that for all $n \geq n_0$,

$$\mu_n(X_\varepsilon(1) \in (0, \lambda_0]) < \left(\frac{\eta}{3C_0}\right)^2.$$

Using $I_{\{S > \varepsilon\}} = I_{\{X_\varepsilon(1) > \lambda_0\}} + I_{\{X_\varepsilon(1) \in (0, \lambda_0]\}}$, $\mathbb{N}_0$-a.e., we have

$$|E_{\mu_n}[F(\vec{X}_{\vec{t}}(\vec{\phi}))I_{\{S>\varepsilon\}}] - E_{\mathbb{N}_0}[F(\vec{X}_{\vec{t}}(\vec{\phi}))I_{\{S>\varepsilon\}}]|$$
(32)
$$\leq |E_{\mu_n}[F(\vec{X}_{\vec{t}}(\vec{\phi}))I_{\{X_\varepsilon(1)>\lambda_0\}}] - E_{\mathbb{N}_0}[F(\vec{X}_{\vec{t}}(\vec{\phi}))I_{\{X_\varepsilon(1)>\lambda_0\}}]|$$
$$+ |E_{\mu_n}[F(\vec{X}_{\vec{t}}(\vec{\phi}))I_{\{X_\varepsilon(1)\in(0,\lambda_0]\}}]|$$
$$+ |E_{\mathbb{N}_0}[F(\vec{X}_{\vec{t}}(\vec{\phi}))I_{\{X_\varepsilon(1)\in(0,\lambda_0]\}}]|.$$

We bound the right-hand side of (32) as follows. By part 2 of Proposition 2.3, the first absolute value is less than $\eta/3$ for $n$ sufficiently large. On the second term, we use the Cauchy–Schwarz inequality to obtain

$$E_{\mu_n}[|F(\vec{X}_{\vec{t}})|I_{\{X_\varepsilon(1)\in(0,\lambda_0]\}}] \leq E_{\mu_n}[F^2(\vec{X}_{\vec{t}}(\vec{1}))]^{1/2} \mu_n(X_\varepsilon(1) \in (0,\lambda_0])^{1/2} \leq \frac{C_0\eta}{3C_0}.$$

The third term is handled similarly. Thus, for $n$ sufficiently large,

$$|E_{\mu_n^\varepsilon}[F(\vec{X}_{\vec{t}}(\vec{\phi}))] - E_{\mathbb{N}_0^\varepsilon}[F(\vec{X}_{\vec{t}}(\vec{\phi}))I_{\{S>\varepsilon\}}]| < \eta,$$

which proves the second condition of Theorem 2.6 for $\{\mu_n^\varepsilon\}_{n\geq 0}$ and $\mathbb{N}_0^\varepsilon$. The result follows by Theorem 2.6. $\square$

**Acknowledgments.** We thank Gordon Slade and Remco van der Hofstad for providing the motivation for this work and for many helpful suggestions. We also thank two anonymous referees for suggestions that led to significant improvements.

EURANDOM  
P.O. Box 513–5600MB  
Eindhoven  
The Netherlands  
E-mail: holmes@eurandom.tue.nl  

Department of Mathematics  
University of British Columbia  
1984 Mathematics Road  
Vancouver, British Columbia  
Canada V6T 1Z2  
E-mail: perkins@math.ubc.ca